\documentclass{amsart}

\usepackage{amsmath, amsthm, amssymb}

\newtheorem{thm}{Theorem}[section]
\newcommand{\bt}{\begin{thm}}
\newcommand{\et}{\end{thm}}

\newtheorem{cor}[thm]{Corollary}
\newcommand{\bc}{\begin{cor}}
\newcommand{\ec}{\end{cor}}

\newtheorem{lem}[thm]{Lemma}
\newcommand{\bl}{\begin{lem}}
\newcommand{\el}{\end{lem}}

\newtheorem{prop}[thm]{Proposition}
\newcommand{\bp}{\begin{prop}}
\newcommand{\ep}{\end{prop}}

\newtheorem{defn}[thm]{Definition}
\newcommand{\bd}{\begin{defn}}      
\newcommand{\ed}{\end{defn}}

\newtheorem{rmrk}[thm]{Remark}
\newcommand{\br}{\begin{rmrk}}
\newcommand{\er}{\end{rmrk}}

\newcommand{\thmref}[1]{Theorem~\ref{#1}}
\newcommand{\secref}[1]{Section~\ref{#1}}
\newcommand{\lemref}[1]{Lemma~\ref{#1}}

\newcommand{\corref}[1]{Corollary~\ref{#1}}
\newcommand{\propref}[1]{Proposition~\ref{#1}}

\newcommand{\N}{\mathbb{N}}

\newcommand{\R}{\mathbb{R}}
\newcommand{\Z}{\mathbb{Z}}
\newcommand{\diam}{\operatorname{diam}}
\newcommand{\length}{\operatorname{length}}
\newcommand{\hm}{{\mathcal H}}

\newcommand{\lip}{\operatorname{Lip}}
\newcommand{\mass}[2][]{{\mathbf M_{#1}}(#2)}

\newcommand{\form}{{\mathcal D}}        
\newcommand{\curr}{{\mathbf M}}         
\newcommand{\norcurr}{{\mathbf N}}      
\newcommand{\rectcurr}{{\mathcal R}}    
\newcommand{\intrectcurr}{{\mathcal I}} 
\newcommand{\intcurr}{{\mathbf I}}      

\newcommand{\rstr}{\:\mbox{\rule{0.1ex}{1.2ex}\rule{1.1ex}{0.1ex}}\:}
\newcommand{\bdry}{\partial\hspace{-0.05cm}}
\newcommand{\slice}[3]{\langle#1,#2,#3\rangle}
\newcommand{\on}[1]{|_{#1}}
\newcommand{\spt}{\operatorname{spt}}
\newcommand{\ohne}{\backslash}

\begin{document}

\title[Isoperimetric inequalities in metric spaces]{Isoperimetric inequalities of euclidean type in metric spaces}

\author{Stefan Wenger}

\address
  {Department of Mathematics,
   ETH Zentrum,
   CH -- 8092 Z\"urich,
   Switzerland}
\email{wenger@math.ethz.ch}

\maketitle

\section{Introduction}
 The purpose of this paper is to prove an isoperimetric inequality of euclidean type for complete 
 metric spaces admitting a cone-type inequality. These include all Banach spaces and all complete, 
 simply-connected metric spaces of non-positive curvature in the sense of Alexandrov or Busemann. 
 The main theorem generalizes results of Gromov \cite{Gromov} and Ambrosio-Kirchheim \cite{Ambr-Kirch-curr}. 

\subsection{Statement of the main result}
The isoperimetric problem of euclidean type for a space $X$ and given classes $\intcurr_{k-1}$, $\intcurr_k$, 
and $\intcurr_{k+1}$ of surfaces of 
dimension $k-1$, $k$, and $k+1$ in $X$, together with boundary operators 
$\intcurr_{k+1}\overset{\partial}{\longrightarrow} \intcurr_k \overset{\partial}{\longrightarrow} \intcurr_{k-1}$ and a volume function
$\mathbf M$ on each class, asks the following: Does there exist for every
surface $T\in \intcurr_k$ without boundary, $\bdry T=0$, a surface $S\in \intcurr_{k+1}$ with $\bdry S=T$ and such that
\begin{equation*}
 \mass{S}\leq D \mass{T}^{\frac{k+1}{k}}
\end{equation*}
for a constant $D$ depending only on $X$ and $k$? A space for which this holds is said to admit an isoperimetric inequality of 
euclidean type for $\intcurr_k$ (or in dimension $k$).
The isoperimetric problem of euclidean type was resolved by Federer and Fleming in \cite{Fed-Flem} 
for euclidean space $X=\R^n$ and in the class $\intcurr_k$ of $k$-dimensional integral currents, $k\in\{1,\dots,n\}$. 
In \cite{Gromov} Gromov extended the result to finite dimensional normed spaces and moreover to complete Riemannian manifolds
admitting a cone-type inequality (for which the definition will be given below). Gromov worked in the class of Lipschitz chains,
formal finite sums of Lipschitz maps on standard simplices.
Recently, Ambrosio and Kirchheim extended in \cite{Ambr-Kirch-curr} the theory of currents from the euclidean setting to general 
metric spaces. The notion of metric integral currents defines suitable classes $\intcurr_k(X)$
of $k$-dimensional surfaces in $X$. (It is to be noted that there are metric spaces for which $\intcurr_k(X)$ only consists of the 
trivial current. However, for the spaces considered here, this is not the case (see below)).
In \cite{Ambr-Kirch-curr} the isoperimetric inequality of euclidean type is proved for dual Banach spaces $X$ admitting an 
approximation by finite dimensional subspaces in the following sense: 
There exists a sequence of projections $P_n: X\longrightarrow X_n$ onto finite dimensional subspaces such that 
$P_n(x)$ weakly$^*$-converges to $x$ for every $x\in X$. 
The authors then raise the question whether all Banach spaces admit an isoperimetric inequality of euclidean type.\\
In this paper we answer this question affirmatively and, in fact, prove the euclidean isoperimetric inequality for a large
class of metric spaces including also many non-linear ones. 
We will work in the class of metric integral currents $\intcurr_k(X)$ developed in \cite{Ambr-Kirch-curr}, the main definitions 
of which will be recalled in \secref{section:currents}. 
As for the moment one can think of a $k$-dimensional integral current as a countable union of bi-Lipschitz surfaces
$f_i : A_i\rightarrow X$ with $A_i\subset \R^k$ compact and with total finite Hausdorff measure.
If the $A_i$ have `regular' boundaries then $\bdry T$ can be viewed as
the union $\bigcup f_i(\bdry A_i)$. An integral current $T$ with $\bdry T= 0$ will be called a cycle.
\bd
 A metric space $(X,d)$ is said to admit a $k$-dimensional cone-type inequality (or to admit a cone-type inequality for $\intcurr_k(X)$)
 if for every cycle $T\in\intcurr_k(X)$ there exists an $S\in\intcurr_{k+1}(X)$ satisfying $\bdry S=T$ and
 \begin{equation*}
  \mass{S}\leq C_k\diam(\spt T)\mass{T}
 \end{equation*}
 for a constant $C_k$ depending only on $k$ and $X$.
\ed
In terms of the intuitive definition given above the support $\spt T$ is the closure of $\bigcup f_i(A_i)$. See 
\secref{section:currents} for precise definitions.
The main result can be stated as follows:
\bt\label{theorem:cone-isoperimetric-inequality}
Let $(X,d)$ be a complete metric space and $k\in\N$. Suppose that $X$ satisfies a cone-type inequality for $\intcurr_k(X)$ 
and,
if $k\geq 2$, that $X$ also satisfies an isoperimetric inequality of euclidean type for $\intcurr_{k-1}(X)$. 
Then $(X,d)$ admits an isoperimetric inequality of euclidean type for $\intcurr_k(X)$: For every cycle $T\in\intcurr_k(X)$ there 
exists an $S\in\intcurr_{k+1}(X)$ with $\bdry S=T$ and such that
\begin{equation*}
 \mass{S}\leq D_k[\mass{T}]^{\frac{k+1}{k}}
\end{equation*}
where $D_k$ only depends on $k$ and the constants of the cone inequality in $\intcurr_k(X)$ and the isoperimetric 
inequality in $\intcurr_{k-1}(X)$.
\et
This theorem extends \cite[3.4.C]{Gromov} from the setting of Riemannian manifolds to that of complete metric spaces.
As will be shown in \secref{subsection:cone-inequality} Banach spaces admit cone-type inequalities in every dimension. This leads
to the following generalization of the result of Ambrosio and Kirchheim and answers raised question.
\bc\label{corollary:banach-isop}
 There are universal constants $D_k$, $k\geq 1$, such that every Banach space $E$ admits an isoperimetric inequality of euclidean type
 for $\intcurr_k(E)$ with constant $D_k$.
\ec
The constants $D_k$ in \thmref{theorem:cone-isoperimetric-inequality} and \corref{corollary:banach-isop} can be computed 
explicitly. However, they are not optimal, not even for $X=\R^n$. The case of optimality for $X=\R^n$ is treated
in \cite{Almgren}.\\
The cone-type inequality is moreover satisfied in spaces $(X,d)$ admitting a $\gamma$-convex bicombing, for some $\gamma>0$. 
By this we mean a choice, for 
all $x,y\in X$, of a rectifiable path $c_{xy}:[0,1]\to X$ joining $x$ to $y$, parameterized proportionally to arc-length and with
$\length(c_{xy})\leq \gamma d(x,y)$, and such that for any three points
$x,y,y'\in X$ the inequality
\begin{equation*}
 d(c_{xy}(t),c_{xy'}(t))\leq t\gamma d(y,y')
\end{equation*}
holds for all $t\in[0,1]$ with a constant $\gamma=\gamma(X)$. 
Examples of such spaces include all simply-connected metric spaces
of non-positive curvature in the sense of Alexandrov and, more general, in the sense of Busemann 
(see \secref{subsection:cone-inequality}). We point out that these spaces contain many rectifiable sets (see \cite{Kleiner}) and hence
$\intcurr_k$ is not trivial.
\bc
 For fixed $\gamma>0$ and $k\in\N$, 
 every complete metric space $(X,d)$ with a $\gamma$-convex bicombing admits an isoperimetric inequality of euclidean type for
 $\intcurr_k(X)$ with constants $D_k$ depending only on $k$ and $\gamma$.
\ec
We only mention one important application of the euclidean isoperimetric inequality here. This is to find 
a solution to the Plateau Problem under suitable 
conditions. \thmref{theorem:cone-isoperimetric-inequality} together with \cite[Theorem 10.6]{Ambr-Kirch-curr} yields the following
theorem.
\bt
If $Y$ is the dual of a separable Banach space then for every $k\in\N$ and each cycle $T\in\intcurr_k(Y)$ with compact support there exists an
$S\in\intcurr_{k+1}(Y)$ with compact support, $\bdry S = T$, and such that
\begin{equation*}
\mass{S}= \inf\{\mass{S'}: \text{$S'\in\intcurr_{k+1}(Y)$, $\bdry S' = T$}\}.
\end{equation*}
\et
This generalizes the corresponding result in \cite{Ambr-Kirch-curr}. However, the Plateau problem remains open in the 
case of general Banach spaces.
\subsection{Outline of the main argument}
The proof of the main theorem is inspired by Gromov's argument. However, the methods in \cite{Gromov} rely in
several ways on the bi-Lipschitz embeddability of compact Riemannian manifolds (and finite dimensional normed spaces) into
euclidean space and do not generalize to the non-embeddable setting. Our approach uses a more intrinsic analysis of
$k$-dimensional cycles.
For the description of our argument it is convenient to introduce the following notation: 
A cycle $T\in\intcurr_k(X)$ is called {\it round} if
\begin{equation*}
  \diam(\spt T)\leq E \mass{T}^{\frac{1}{k}}
\end{equation*}
for a constant $E$ depending only on $k$ and on the space $X$. 
The essential step in the proof is stated in \propref{proposition:effective-decomposition} which claims the following:
Under the hypotheses of \thmref{theorem:cone-isoperimetric-inequality} every cycle $T\in\intcurr_k(X)$ can be decomposed into the sum $T=T_1+\dots+T_N+R$ of round cycles $T_i$ and a cycle $R$ 
with the properties that $\sum\mass{T_i}\leq(1+\lambda)\mass{T}$ and $\mass{R}\leq(1-\delta)\mass{T}$ for constants $0<\delta,\lambda<1$
depending only on $k$ and the constant from the isoperimetric inequality for $\intcurr_{k-1}(X)$.
The construction of such a decomposition is based on an analysis of the growth of the function 
$r\mapsto \|T\|(B(y,r))$ .
In the intuitive language used above, $\|T\|(B(y,r))$ is the volume of the intersection of $\bigcup f_i(A_i)$ with the ball with 
radius $r$ and center $y$. For almost every $y\in \bigcup f_i(A_i)$ we have
\begin{equation}\label{equation:growth-intro}
 \|T\|(B(y,r))\geq Fr^k
\end{equation}
for small $r>0$ and a constant $F$ depending only on $k$. 
Denoting by $r_0(y)$ the least upper bound of those $r$ satisfying \eqref{equation:growth-intro} 
one can prove that, when $T$ is cut
open along the metric sphere with center $y$ and of radius about $r_0(y)$ only little boundary is created. By closing this 
boundary with a suitable isoperimetric filling (\lemref{lemma:isoperimetric-support}) one constructs a decomposition of $T$ into a
sum $T=T_1+\tilde{R}$. The cycle $T_1$, lying essentially in $B(y,r_0(y))$ is round. This will easily follow from the definition
of $r_0(y)$. By using a simple Vitali-type covering argument one then shows that enough such round cycles $T_i$ 
can be split off in order to leave a rest $R$ satisfying $\mass{R}\leq (1-\delta)\mass{T}$.
Successive application of \propref{proposition:effective-decomposition} will easily establish the proof of the main theorem. 

\medskip

The paper is structured as follows: In \secref{section:currents} we recall the main definitions from the theory of 
currents in metric spaces and state those results from \cite{Ambr-Kirch-curr} vital for our purposes. Then, following a construction
in \cite{Ambr-Kirch-curr}, we prove cone-type inequalities for various classes of metric spaces.
The decomposition of a given cycle is constructed in \secref{section:partial-decomposition}. This forms the main part 
of the paper. \secref{section:partial-decomposition}  also contains the proof of the main theorem.\\
\section{Currents in metric spaces and cone constructions}\label{section:currents}
This section contains the main definitions from the theory of metric currents developed in \cite{Ambr-Kirch-curr} as well
as the results relevant in our context. The purpose of \secref{subsection:cone-inequality} is to construct cone fillings
and to prove a cone-type inequality for metric spaces admitting a $\gamma$-convex bicombing. The latter condition is reviewed
and a list of examples is given.
\subsection{Definitions and theorems}
{\sloppy
Let $(X,d)$ be a complete metric space and let $\form^k(X)$ denote the set of $(k+1)$-tuples $(f,\pi_1,\dots,\pi_k)$ 
of Lipschitz functions on $X$ with $f$ bounded. The Lipschitz constant of a Lipschitz function $f$ on $X$ will
be denoted by $\lip(f)$.
}
\bd
A $k$-dimensional metric current  $T$ on $X$ is a multi-linear functional on $\form^k(X)$ satisfying the following
properties:
\begin{enumerate}
 \item If $\pi^j_i$ converges point-wise to $\pi_i$ as $j\to\infty$ and if $\sup_{i,j}\lip(\pi^j_i)<\infty$ then
       \begin{equation*}
         T(f,\pi^j_1,\dots,\pi^j_k) \longrightarrow T(f,\pi_1,\dots,\pi_k).
       \end{equation*}
 \item If $\{x\in X:f(x)\not=0\}$ is contained in the union $\bigcup_{i=1}^kB_i$ of Borel sets $B_i$ and if $\pi_i$ is constant 
       on $B_i$ then
       \begin{equation*}
         T(f,\pi_1,\dots,\pi_k)=0.
       \end{equation*}
 \item There exists a finite Borel measure $\mu$ on $X$ such that
       \begin{equation}\label{equation:mass-def}
        |T(f,\pi_1,\dots,\pi_k)|\leq \prod_{i=1}^k\lip(\pi_i)\int_X|f|d\mu
       \end{equation}
       for all $(f,\pi_1,\dots,\pi_k)\in\form^k(X)$.
\end{enumerate}
\ed
The space of $k$-dimensional metric currents on $X$ is denoted by $\curr_k(X)$ and the minimal Borel measure $\mu$
satisfying \eqref{equation:mass-def} is called mass of $T$ and written as $\|T\|$. 
The support of $T$ is, by definition, the closed set $\spt T$ of points $x\in X$ such that $\|T\|(B(x,r))>0$ for all $r>0$. 
Hereby, $B(x,r)$ denotes the closed ball $B(x,r):= \{y\in X : d(y,x)\leq r\}$.
\br
As is done in \cite{Ambr-Kirch-curr} we will also assume here
that $\spt T$ is separable and furthermore that $\|T\|$ is concentrated on a $\sigma$-compact set, i.\ e.\ $\|T\|(X\ohne C) = 0$
for a $\sigma$-compact set $C\subset X$.
However, as is observed in \cite{Ambr-Kirch-curr}, these facts can be proved if one accepts the standard ZFC-set theory.
\er
%
%
The restriction of $T\in\curr_k(X)$ to a Borel set $A\subset X$ is given by 
\begin{equation*}
  (T\rstr A)(f,\pi_1,\dots,\pi_k):= T(f\chi_A,\pi_1,\dots,\pi_k).
\end{equation*}
This expression is well-defined since $T$ can be extended to a functional on tuples for which the first argument lies in 
$L^\infty(X,\|T\|)$.\\
%
The boundary of $T\in\curr_k(X)$ is the functional
\begin{equation*}
 \bdry T(f,\pi_1,\dots,\pi_{k-1}):= T(1,f,\pi_1,\dots,\pi_{k-1}).
\end{equation*}
It is clear that $\bdry T$ satisfies conditions (i) and (ii) in the above definition. If $\bdry T$ also has 
finite mass (condition (iii)) then $T$ is called a normal current. The respective space is denoted by $\norcurr_k(X)$.\\
The push-forward of $T\in\curr_k(X)$ 
under a Lipschitz map $\varphi$ from $X$ to another complete metric space $Y$ is given by
\begin{equation*}
 \varphi_\# T(g,\tau_1,\dots,\tau_k):= T(g\circ\varphi, \tau_1\circ\varphi,\dots,\tau_k\circ\varphi)
\end{equation*}
for $(g,\tau_1,\dots,\tau_k)\in\form^k(Y)$. This defines a $k$-dimensional current on $Y$.\\
In this paper we will mainly be concerned with integer rectifiable and integral currents.
For notational purposes we first repeat some well-known definitions. The Hausdorff $k$-dimensional 
measure of $A\subset X$ is defined to be
\begin{equation*}
 \hm^k(A):= \lim_{\delta\searrow 0}\inf\left\{\sum_{i=1}^\infty \omega_k\left(\frac{\diam(B_i)}{2}\right)^k :
      B\subset \bigcup_{i=1}^\infty B_i\text{, }\diam(B_i)<\delta\right\},
\end{equation*}
where $\omega_k$ denotes the Lebesgue measure of the unit ball in $\R^k$. The $k$-dimensional 
lower density $\Theta_{*k}(\mu, x)$ of a finite Borel measure $\mu$ at a point $x$ is given by the formula
\begin{equation*}
 \Theta_{*k}(\mu, x):= \liminf_{r\searrow 0}\frac{\mu(B(x,r))}{\omega_k r^k}.
\end{equation*}
An $\hm^k$-measurable set $A\subset X$
is said to be countably $\hm^k$-rectifiable if there exist Lipschitz maps $f_i :B_i\longrightarrow X$ from subsets
$B_i\subset \R^k$ such that
\begin{equation*}
\hm^k(A\ohne \bigcup f_i(B_i))=0.
\end{equation*}
\bd
A current $T\in\curr_k(X)$ with $k\geq 1$ is said to be rectifiable if
\begin{enumerate}
 \item $\|T\|$ is concentrated on a countably $\hm^k$-rectifiable set and
 \item $\|T\|$ vanishes on $\hm^k$-negligible sets.
\end{enumerate}
$T$ is called integer rectifiable if, in addition, the following property holds:
\begin{enumerate}
 \item[(iii)] For any Lipschitz map $\varphi\colon X\longrightarrow \R^k$ and any open set $U\subset X$ there exists 
       $\theta\in L^1(\R^k,\Z)$ such that
       \begin{equation*}
        \varphi_\#(T\rstr U)(f,\pi_1,\dots,\pi_k)= \int_{\R^k}\theta f\det\left(\frac{\partial\pi_i}{\partial x_j}\right)d{\mathcal L}^k
       \end{equation*}
       for all $(f,\pi_1,\dots,\pi_k)\in\form^k(\R^k)$.
\end{enumerate}
\ed
A $0$-dimensional (integer) rectifiable current is a $T\in\curr_0(X)$ of the form
\begin{equation*}
 T(f)=\sum_{i=1}^\infty \theta_i f(x_i),\qquad \text{$f$ Lipschitz and bounded,}
\end{equation*}
for suitable $\theta_i\in\R$ (or $\theta_i\in \Z$) and $x_i\in X$.\\
The space of rectifiable currents is denoted by $\rectcurr_k(X)$, that of integer rectifiable currents by $\intrectcurr_k(X)$.
As is easily seen, $\intrectcurr_k(X)$ is a Banach subspace of $\curr_k(X)$ endowed with the mass norm.
Integer rectifiable normal currents are called integral currents. The respective space is denoted by $\intcurr_k(X)$.
In the following, an element $T\in\intcurr_k(X)$ with zero boundary $\bdry T=0$ will be called a cycle.\\
%
%
The characteristic set $S_T$ of a rectifiable current $T\in\rectcurr_k(X)$ is defined by
\begin{equation}\label{equation:characteristic-set}
 S_T:= \{x\in X: \Theta_{\star k}(\|T\|, x)>0\}.
\end{equation}
It can be shown that $S_T$ is countably $\hm^k$-rectifiable and that $\|T\|$ is concentrated on $S_T$.
\begin{thm}[{\cite[Theorem 9.5]{Ambr-Kirch-curr}}]\label{theorem:mass-representation}
 If $T\in\rectcurr_k(X)$ then there exist $\hm^k$-integrable functions $\lambda: S_T\longrightarrow [k^{-k/2},2^k/\omega_k]$ and 
 $\theta: S_T\longrightarrow (0,\infty)$ such that
 \begin{equation*}
  \|T\|(A)=\int_{A\cap S_T}\lambda\theta d\hm^k\qquad\text{for $A\subset X$ Borel,}
 \end{equation*}
 that is, $\|T\|= \lambda\theta d\hm^k\rstr S_T$.
 Moreover, if $T$ is an integral current then $\theta$ takes values in $\N:=\{1,2,\dots\}$ only.
\end{thm}
The following Slicing Theorem (proved in \cite[Theorems 5.6 and 5.7]{Ambr-Kirch-curr})
is, beside \thmref{theorem:mass-representation}, the only result from the theory of
metric currents needed in the proof of the main result.
\bt\label{theorem:slicing}
Let be $T\in\norcurr_k(X)$ and $\varrho$ a Lipschitz function on $X$. Then there exists for almost every $r\in\R$ a normal current
$\slice{T}{\varrho}{r}\in\norcurr_{k-1}(X)$ with the following properties:
\begin{enumerate}
 \item $\slice{T}{\varrho}{r}= \bdry(T\rstr\{\varrho \leq r\}) - (\bdry T)\rstr\{\varrho\leq r\}$
 \item $\|\slice{T}{\varrho}{r}\|$ and $\|\bdry\slice{T}{\varrho}{r}\|$ are concentrated on $\varrho^{-1}(\{r\})$
 \item $\mass{\slice{T}{\varrho}{r}}\leq\lip(\varrho)\frac{d}{dr}\mass{T\rstr\{\varrho\leq r\}}$.
\end{enumerate}
Moreover, if $T\in\intcurr_k(X)$ then $\slice{T}{\varrho}{r}\in\intcurr_{k-1}(X)$ for almost all $r\in\R$.
\et
\subsection{Cone-type inequalities}\label{subsection:cone-inequality}
The following cone construction is a slightly modified version of the one given in \cite{Ambr-Kirch-curr}.\\
Let $(X,d)$ be a complete metric space and $T\in\norcurr_k(X)$ and endow $[0,1]\times X$ with the euclidean product metric.
Given a Lipschitz function $f$ on $[0,1]\times X$ and $t\in[0,1]$ we define the function $f_t:X\longrightarrow \R$ by
$f_t(X):= f(t,x)$. To every $T\in\norcurr_k([0,1]\times X)$ and every $t\in[0,1]$ we associate the normal $k$-current on $X$ given
by the formula
\begin{equation*}
  ([t]\times T)(f,\pi_1,\dots,\pi_k):= T(f_{t},\pi_{1\,t},\dots,\pi_{k\,t}),
\end{equation*}
The product of a normal current with the interval $[0,1]$ is defined 
as follows.
 \bd
  For a normal current $T\in\norcurr_k(X)$
  the functional $[0,1]\times T$ on $\form^{k+1}([0,1]\times X)$ is given by
  \begin{equation*}
  \begin{split}
   [0,1]\times& T (f,\pi_1,\dots,\pi_{k+1}):= \\
      &\sum_{i=1}^{k+1}(-1)^{i+1}\int_0^1T\left(f_t\frac{\partial \pi_{i\,t}}{\partial t},\pi_{1\,t},
                                               \dots,\pi_{i-1\,t},\pi_{i+1\,t},\dots,\pi_{k+1\,t}\right)dt
  \end{split}
  \end{equation*}
  for $(f,\pi_1,\dots,\pi_{k+1})\in\form^{k+1}([0,1]\times X)$.
 \ed
It is easily seen that this definition is well-posed. Moreover, we have the following result whose proof is analogous to that
of \cite[Proposition 10.2 and Theorem 10.4]{Ambr-Kirch-curr}.
\bt\label{theorem:cone-construction}
  For every $T\in\norcurr_k(X)$ with bounded support the functional $[0,1]\times T$ is a $(k+1)$-dimensional normal current
 on $[0,1]\times X$ with boundary
  \begin{equation*}
   \partial([0,1]\times T)= [1]\times T - [0]\times T + [0,1]\times\partial T.
  \end{equation*}
  Moreover, if $T\in\intcurr_k(X)$ then $[0,1]\times T\in\intcurr_{k+1}([0,1]\times X)$.
\et
We now define $\gamma$-convex bicombings and give a list of examples of spaces sharing this property. 
Recall that a rectifiable curve is a continuous map $c: [0,1]\longrightarrow X$ with finite length, i.\ e.\ such that
 \begin{equation*}
  \length(c):= \sup\left\{\sum_{i=1}^N d(c(t_i),c(t_{i+1})) : 0=t_1<\dots<t_N=1\right\}<\infty.
 \end{equation*}
 \bd
  A $\gamma$-convex bicombing on $X$ is a choice, for each two points $x,y\in X$, of a rectifiable
  curve $c_{xy}: [0,1]\to X$ joining $x$ to $y$ such that the following conditions hold:
  \begin{enumerate}
   \item $c_{xy}$ is parameterized proportional to arc-length:
    \begin{equation*}
     \length(c_{xy}\on{[0,s]})= s\length(c_{xy})\quad\text{for all $s\in[0,1]$}.
    \end{equation*}
   \item $\length(c_{xy})\leq \gamma d(x,y)$.
   \item For any three points $x,y,y'\in X$ and for $t\in[0,1]$ we have
  \begin{equation*}
   d(c_{xy}(t), c_{xy'}(t))\leq t\gamma d(y,y').
  \end{equation*}
 \end{enumerate}
 \ed
 The following spaces admit a $1$-convex bicombings:\\
 1. Every Banach space together with $c_{xy}(t):= ty+(1-t)x$.\\
 2. Every Hadamard manifold (i.\ e.\ complete, simply-connected, with non-positive sectional curvature) 
         with $c_{xy}:=$ the unique geodesic (parameterized on $[0,1]$) from $x$ to $y$.\\
 3. More generally, 2.\ holds for every simply-connected complete metric space of non-positive curvature in 
         the sense of Alexandrov or (even more generally) in the sense of Busemann.
    For an account on these spaces see \cite{Bridson-Haefliger}.\\
 4. The space $H(E)$ of bounded, closed, convex subsets of a given Banach space $E$ endowed with the Hausdorff metric. 
    A convex bicombing is given by $c_{AA'}(t):= (1-t)A+tA'$, $A,A'\in H(E)$.
    \\
 \bp
  If $(X,d)$ is a complete metric space admitting a $\gamma$-convex bicombing then every cycle $T\in\intcurr_k(X)$ has a filling 
  $S\in\intcurr_{k+1}(X)$ satisfying
  \begin{equation*}
   \mass{S}\leq (k+1)\gamma^{k+1}\diam(\spt T)\mass{T}.
  \end{equation*}
 \ep
 \begin{proof}
  We fix $x_0\in\spt T$ and define a locally Lipschitz map $\varphi:[0,1]\times X\to X$ by $\varphi(t,x):= c_{x_0x}(t)$. It is clear that 
  for fixed $x\in\spt T$ the map $t\mapsto \varphi(t,x)$ is $\gamma\diam(\spt T)$-Lipschitz, whereas for fixed $t\in[0,1]$ 
  the map $x\mapsto \varphi(t,x)$ is $\gamma$-Lipschitz. \thmref{theorem:cone-construction} implies that
  $\varphi_{\#}([0,1]\times T)\in\intcurr_{k+1}(X)$ and furthermore
  \begin{equation*}
    \bdry \varphi_{\#}([0,1]\times T) = \varphi_{\#}(\bdry([0,1]\times T))
                                      = \varphi_{\#}([1]\times T) - \varphi_{\#}([0]\times T)
                                      = T.
  \end{equation*}
  To obtain the estimate on mass we compute for $(f,\pi_1,\dots,\pi_{k+1})\in\form^{k+1}(X)$ that
  \begin{equation*}
   \begin{split}
    |\varphi_{\#}(&[0,1]\times T)(f,\pi_1,\dots,\pi_{k+1})| \\
       &\leq \sum_{i=1}^{k+1}\left|\int_0^1 
             T(f\circ\varphi_t\frac{\partial \pi_i\circ\varphi_t}{\partial t},\pi_1\circ\varphi_t,\dots,\pi_{i-1}\circ\varphi_t,
               \pi_{i+1}\circ\varphi_t,\dots,\pi_{k+1}\circ\varphi_t)dt\right|\\
                                        &\leq \sum_{i=1}^{k+1} \int_0^1\prod_{j\not=i}\lip(\pi_j\circ\varphi_t)
                                                \int_X\left|f\circ\varphi_t\frac{\partial (\pi_i\circ\varphi_t)}{\partial t}\right|
                                                 d\|T\|dt\\
                                        &\leq (k+1)\gamma^{k+1}\diam(\spt T)\prod_{j=1}^{k+1}\lip(\pi_j)\int_0^1
                                                \int_X|f\circ\varphi(t,x)|d\|T\|(x)dt.
   \end{split}
  \end{equation*}
  From this it follows that $\|[0,1]\times T\|\leq (k+1)\gamma^{k+1}\diam(\spt T)\varphi_{\#}({\mathcal L}^1\times\|T\|)$ and this concludes
  the proof.
 \end{proof} 
\section{Partial decomposition and proof of the main result}\label{section:partial-decomposition}
The aim of this section is to prove the proposition below which forms the crucial step when decomposing a cycle into the sum
of round cycles. This result will be used to prove the \thmref{theorem:cone-isoperimetric-inequality}.
\bp\label{proposition:effective-decomposition}
 Let $(X,d)$ be a complete metric space and $k\geq 1$ an integer. If $k\geq 2$ then suppose furthermore that $X$ has a euclidean
isoperimetric inequality for $\intcurr_{k-1}(X)$ with a constant $C>0$. There then exist constants 
$E>0$ and $0<\delta,\lambda<1$ depending only on $k$ and $C$ with the following property:
Every cycle $T\in\intcurr_k(X)$ admits a decomposition  $T= \sum_{i=1}^N T_i + R$ into a sum of integral 
 cycles satisfying:
 \begin{enumerate}
  \item $\diam(\spt T_i)\leq E\mass{T_i}^{1/k}$
  \item $\mass{R}\leq (1-\delta)\mass{T}$
  \item $\sum_{i=1}^N\mass{T_i}\leq (1+\lambda)\mass{T}$.
 \end{enumerate}
\ep
We first state some preparatory lemmas. The first will be employed to obtain the estimate in 
(ii) for the cycle $R$.
As it is an analog to the simple Vitali covering lemma the proof is omitted.
 \bl\label{lemma:partial-cover}
  Let $(Y,d)$ be a metric space, $\mu$ a finite Borel measure on $Y$, and $F>0$, $k\in\N$. For $y\in Y$ define
  \begin{equation*}
   r_0(y):= \max\{r\geq 0 : \mu(B(y,r))\geq F r^k\}.
  \end{equation*}
  If $r_0(y)>0$ for $\mu$-almost every $y\in Y$  then there exist points $y_1,\dots, y_N\in Y$ satisfying
  \begin{enumerate}
   \item $r_0(y_i)>0$
   \item $B\left(y_i, 2r_0(y_i)\right)\cap B\left(y_j, 2r_0(y_j)\right)=\emptyset$ if $i\not= j$
   \item $\sum_{i=1}^N\mu(B(y_i,r_0(y_i)))\geq \alpha\mu(Y)$
  \end{enumerate}
  for a constant $\alpha>0$ depending only on $k$.
 \el
%
%
The study of the growth of the function $r\mapsto \|T\|(B(x,r))$ will play a predominant role in the proof of 
\propref{proposition:effective-decomposition}. In this context the following easy fact will be helpful.

 \bl\label{lemma:polynomial}
  Fix $\bar{C}>0$, $k\geq 2$, $0\leq r_0<r_1<\infty$, and suppose $\beta : [r_0,r_1]\longrightarrow (0,\infty)$ is non-decreasing and satisfies
  \begin{enumerate}
   \item $\beta(r_0)= \frac{r_0^k}{\bar{C}^{k-1}k^k}$
   \item $\beta(r)\leq \bar{C}[\beta'(r)]^{k/(k-1)}$ for a.\ e.\ $r\in(r_0,r_1)$.
  \end{enumerate}
  Then it follows that
  \begin{equation*}
   \beta(r)\geq \frac{r^k}{\bar{C}^{k-1}k^k}\quad\text{ for all } r\in[r_0,r_1].
  \end{equation*}
 \el
\begin{proof}
 By rearranging (ii) we obtain
 \begin{equation*}
  \frac{\beta'(t)}{\beta(t)^{\frac{k-1}{k}}}\geq \frac{1}{\bar{C}^{\frac{k-1}{k}}}
 \end{equation*}
 and integration from $r_0$ to $r$ yields the claimed estimate.
%
 \end{proof}
The next statement is concerned with the support of isoperimetric fillings. It will be used to prove the roundness of the cycles
$T_i$. The $T_i$ will be constructed by restricting $T$ to a ball $B(y_i,r)$ and filling in the boundary $\bdry(T\rstr B(y_i,r))$ by an
isoperimetric filling. \lemref{lemma:isoperimetric-support} ensures that we can choose a filling whose support stays near its boundary.
\bl\label{lemma:isoperimetric-support}
Let $(X,d)$ be a complete metric space and $k\geq 2$. Suppose that $X$ admits a euclidean isoperimetric inequality for 
$\intcurr_{k-1}(X)$ with a constant $C>0$.
Then there exists for every cycle $T\in\intcurr_{k-1}(X)$ an isoperimetric filling $S\in\intcurr_k(X)$ of $T$ satisfying
\begin{equation*}
 \spt S\subset B(\spt T, 3Ck\mass{T}^{\frac{1}{k-1}}).
\end{equation*}
\el
Hereby, $B(A, \varrho)$ denotes the $\varrho$-neighborhood of the set $A$.
The proof of the lemma is essentially contained in the proof of \cite[Theorem 10.6]{Ambr-Kirch-curr}.
\begin{proof}
Let ${\mathcal M}$ denote the complete metric space consisting of all fillings $S\in\intcurr_k(X)$ of $T$ and endowed with the 
metric given by $d_{\mathcal M}(S,S'):= \mass{S-S'}$. By the Ekeland-Bishop-Phelps variational principle there exists an 
$S\in{\mathcal M}$ satisfying the isoperimetric inequality and such that
the function
\begin{equation*}
 S'\mapsto \mass{S'}+\frac{1}{2}\mass{S'-S}
\end{equation*}
is minimal at $S'=S$. Let be $x\in\spt S\ohne\spt T$ and set $\varrho_x(y):= d(x,y)$. Then, for almost every $0<r<d(x,\spt T)$
the slice $\slice{S}{\varrho_x}{r}$ exists, has zero boundary, and belongs to $\intcurr_{k-1}(X)$. For an
isoperimetric filling $S_r\in\intcurr_k(X)$ of $\slice{S}{\varrho_x}{r}$ the integral current
$S\rstr B^c(x,r)+S_r$ has boundary $T$ and thus, comparison with $S$ yields
\begin{equation*}
\mass{S\rstr B^c(x,r)+S_r} + \frac{1}{2}\mass{S\rstr B(x,r)-S_r}\geq \mass{S}.
\end{equation*}
Hereby, $B^c(x,r)$ denotes the complement of the ball $B(x,r)$.
Together with the isoperimetric inequality, the above estimate imples that
\begin{equation*}
\mass{S\rstr B(x,r)}\leq 3\mass{S_r}\leq 3C\mass{\slice{S}{\varrho_x}{r}}^{\frac{k}{k-1}}
                                                   \text{ for a.\ e.\ $r\in(0,d(x,\spt T))$}.
\end{equation*}
Setting $\beta(r):= \|S\|(B(x,r))$ and using the Slicing Theorem we obtain the inequality
\begin{equation*}
 \beta(r)\leq 3C[\beta'(r)]^{\frac{k}{k-1}}
    \text{ for a.\ e.\ $r\in(0,d(x,\spt T))$}
\end{equation*}
which, after applying \lemref{lemma:polynomial}, yields the claimed bound on $d(x,\spt T)$.
\end{proof}
We are now ready to prove the proposition.
\begin{proof}[Proof of \propref{proposition:effective-decomposition}]
First of all, let $F$ be given by $F:= \frac{\lambda^{k-1}}{C^{k-1}k^k}$ with $\lambda\leq\frac{1}{6}$ small enough such that
$F<\frac{\omega_k}{k^{\frac{k}{2}}}$. Let $T\in\intcurr_k(X)$ be a cycle and define for $y\in X$
\begin{equation*}
r_0(y):= \max\left\{r\geq 0: \|T\|(B(y,r))\geq Fr^k\right\}.
\end{equation*}
By \cite[Theorem 9]{Kirchheim} we have
\begin{equation*}
 \lim_{r\searrow 0}\frac{\hm^k(B(y,r))}{\omega_k r^k} = 1
\end{equation*}
for $\hm^k$-almost all $y\in S_T$, the set $S_T$ being defined as in \eqref{equation:characteristic-set}. Together with 
\thmref{theorem:mass-representation} this implies that the set $Y$ of points $y\in S_T$
satisfying $r_0(y)>0$ has full $\|T\|$-measure.
By \lemref{lemma:partial-cover} there exist points $y_1,\dots, y_N\in Y$ with  $r_0(y_i)>0$, and
such that the balls $B(y_i, 2r_0(y_i))$ are pairwise disjoint and satisfy
\begin{equation}\label{eq:vitali-cover-estimate}
\sum_{i=1}^N\|T\|(B(y_i, r_0(y_i)))\geq \alpha \|T\|(Y) = \alpha\mass{T}
\end{equation}
for a constant $\alpha>0$ depending only on $k$.
Fix $i\in\{1,\dots,N\}$ and set $r_0:=r_0(y_i)$ and $\beta(r):= \|T\|(B(y_i,r))$. It is clear that $\beta$ is non-decreasing,
that $\beta(r_0)=Fr_0^k$, and that $\beta(r)<Fr^k$ for all $r>r_0$. Denote furthermore by $\varrho$ the function 
$\varrho(x):= d(y_i,x)$. By \thmref{theorem:slicing} the slice $\slice{T}{\varrho}{r}=\bdry(T\rstr B(y_i,r))$ exists for almost all $r$, is an element
of $\intcurr_{k-1}(X)$, and satisfies moreover
\begin{equation}\label{equation:slice-derivative}
\mass{\slice{T}{\varrho}{r}}\leq \beta'(r)\quad\text{for a.e. $r$}.
\end{equation}
We now consider one dimensional and higher dimensional cycles separately: If $k=1$ it follows from 
the fact that $F=1$ and from the definition
of $\beta$ that there exists a measurable set $\Omega\subset[r_0,2r_0)$ of positive measure and such that 
$\beta'(r)<1$ for all $r\in\Omega$. Since, for $r\in\Omega$, the slice $\slice{T}{\varrho}{r}$ is a $0$-dimensional integral current, 
$\mass{\slice{T}{\varrho}{r}}$
is an integer number and hence, by \eqref{equation:slice-derivative}, the integral current $T_i:= T\rstr B(y_i, r)$ 
has zero boundary.
Applying this to each $i\in\{1,\dots,N\}$ one easily obtains
a decomposition $T=T_1+\dots+ T_N+R$
satisfying all the properties stated in the proposition (with $\lambda=0$, $\delta\geq\alpha$, and $E<4$).\\
If $k\geq 2$ then \lemref{lemma:polynomial} and the definitions of $F$ and $r_0$ imply the existence of $\Omega\subset[r_0,\frac{4r_0}{3}]$ 
of positive measure such that
\begin{equation}\label{eq:beta-growth}
 C[\beta'(r)]^{\frac{k}{k-1}}<\lambda\beta(r)\quad\text{for all $r\in\Omega$}.
\end{equation}
By \thmref{theorem:slicing} we can assume without loss of generality that the slice $\slice{T}{\varrho}{r}$ exists for every
$r\in\Omega$ and are elements of $\intcurr_{k-1}(X)$. Choose an $r\in\Omega$ arbitrarily and a filling $S\in\intcurr_k(X)$ 
of $\slice{T}{\varrho}{r}$ as in \lemref{lemma:isoperimetric-support}. Together with \eqref{equation:slice-derivative} and 
\eqref{eq:beta-growth} this implies that
\begin{equation}\label{eq:fill-estimate}
\mass{S}\leq\lambda\beta(r)
\end{equation}
and furthermore, since $\lambda\leq \frac{1}{6}$, that the support of $S$ lies in the ball with center $y_i$ and 
with a radius $\bar{r}$ satisfying
\begin{equation*}
\bar{r}\leq \frac{4}{3}r_0+3Ck[\mass{\slice{T}{\varrho}{r}}]^\frac{1}{k-1}
        \leq \frac{4}{3}\left(1 + \frac{3Ck(\lambda F)^{\frac{1}{k}}}{C^{\frac{1}{k}}}\right)r_0
       \leq 2r_0.
\end{equation*}
Clearly, $T_i:= T\rstr B(y_i,r) - S$ defines an integral cycle which satisfies
\begin{equation*}
 (1-\lambda)\beta(r)\leq\mass{T_i}\leq(1+\lambda)\beta(r).
\end{equation*}
Since $T_i$ has support in $B(y_i, 2r_0(y_i))$ it follows that
\begin{equation*}
\diam(\spt T_i)\leq 4r_0(y_i) =  \frac{4}{F^{\frac{1}{k}}}[\beta(r_0)]^\frac{1}{k} 
               \leq \frac{4}{[F(1-\lambda)]^\frac{1}{k}}\mass{T_i}^\frac{1}{k}
\end{equation*}
and hence $T_i$ fulfills condition (i) with $E:= \frac{4}{[F(1-\lambda)]^\frac{1}{k}}$.\\
Since our construction of $T_i$ leaves $T\rstr B^c(y_i, 2r_0(y_i))$ unaffected (by the fact that the balls $B(y_j,2r_0(y_j))$ are 
pairwise disjoint) we can apply the above construction to every $i\in\{1,\dots, N\}$ to obtain round cycles $T_1,\dots,T_N$. Setting
$R:= T- \sum_{i=1}^N T_i$ this yields a decomposition $T=T_1+\dots +T_N+R$ satisfying the claimed properties.
Indeed, we have
\begin{equation*}
 \sum_{i=1}^N\mass{T_i}\leq (1+\lambda)\sum_{i=1}^N \|T\|(B_i)\leq (1+\lambda)\mass{T}
\end{equation*}
where $B_i$ is the ball chosen individually for every $i$ as above. The estimate for $\mass{R}$ is also obvious since,
by \eqref{eq:vitali-cover-estimate} and \eqref{eq:fill-estimate}, we have
\begin{equation*}
 \mass{R}\leq \|T\|(X\ohne \bigcup B_i) + \lambda \sum\|T\|(B_i)
         \leq (1-\alpha(1-\lambda))\mass{T}.
\end{equation*}
This completes the proof of the proposition with $\delta:= \alpha(1-\lambda)$.
\end{proof}
The isoperimetric inequality now easily follows from \propref{proposition:effective-decomposition}.
\begin{proof}[Proof of \thmref{theorem:cone-isoperimetric-inequality}]
 Let $T\in\intcurr_k(X)$ be a cycle. 
 Successive application of \propref{proposition:effective-decomposition} yields (possibly finite)
 sequences of cycles $(T_i)_i$, $(R_n)_n\subset\intcurr_k(X)$ and an increasing sequence $(N_n)_n\subset\N$
 with the following properties:
 \begin{itemize}
  \item $T= \sum_{i=1}^{N_n}T_i + R_n$
  \item $\diam(\spt T_i)\leq E\mass{T_i}^{1/k}$
  \item $\mass{R_n}\leq (1-\delta)^n\mass{T}$
  \item $\sum_{i=1}^\infty \mass{T_i}\leq \left[(1+\lambda)\sum_{i=0}^\infty(1-\delta)^i\right] \mass{T}=\frac{1+\lambda}{\delta}\mass{T}$.
 \end{itemize}

 The isoperimetric filling of $T$ is then constructed as follows. We first fill each
 $T_i$ with an $S_i\in\intcurr_{k+1}(X)$ from the cone inequality, i.\ e.\ one with $\bdry S_i = T_i$ and such that
 \begin{equation}\label{eq:mass-S}
  \mass{S_i}\leq C_k\diam(\spt T_i)\mass{T_i} \leq C_kE\mass{T_i}^{\frac{k+1}{k}}.
 \end{equation}
 The finiteness of $\sum_{i=1}^\infty\mass{T_i}$ implies that the sequence $S^n:= \sum_{i=1}^{N_n}S_i$
 is a Cauchy-sequence with respect to the mass norm because
 \begin{equation*}
  \mass{S^{n+q}-S^n}\leq C_kE\sum_{i=N_n+1}^\infty\mass{T_i}^{\frac{k+1}{k}}
                    \leq C_kE\left[\sum_{i=N_n+1}^\infty\mass{T_i}\right]^{\frac{k+1}{k}}
\;\;\overset{n\rightarrow\infty}{\longrightarrow}0.
 \end{equation*}
 Since $\intrectcurr_{k+1}(X)$ is a Banach space the sequence $S^n\in\intcurr_{k+1}(X)\subset\intrectcurr_{k+1}(X)$ 
 converges to a limit current $S\in\intrectcurr_{k+1}(X)$. As $T-\bdry S^n=R_n$ converges to $0$ it follows that
 $\bdry S=T$ and, in particular, that $S\in\intcurr_{k+1}(X)$.
 Finally, $S$ is an isoperimetric filling of $T$. Indeed, we have
\begin{equation*}
 \mass{S}\leq \sum\mass{S_i}\leq C_kE\sum\mass{T_i}^{\frac{k+1}{k}}
         \leq C_kE\left(\frac{1+\lambda}{\delta}\right)^\frac{k+1}{k}\mass{T}^{\frac{k+1}{k}},
\end{equation*}
which completes the proof.
\end{proof}
%
%
%
 
%
\end{document}